\documentclass[12pt,a4paper,twoside]{article}

\usepackage{amsmath, amssymb}
\usepackage{ascmac}

\DeclareMathOperator{\Pic}{Pic}

\DeclareMathOperator{\Cliff}{Cliff}
\DeclareMathOperator{\gon}{gon}

\DeclareMathOperator{\End}{End}
\DeclareMathOperator{\Hom}{Hom}
\DeclareMathOperator{\Grass}{Grass}
\DeclareMathOperator{\rk}{rk}
\begin{document}

\title{\textbf{\Large{Lifts of line bundles on curves on K3 surfaces}}}

\author{Kenta Watanabe and Jiryo Komeda \thanks{Nihon University, College of Science and Technology,   7-24-1 Narashinodai Funabashi city Chiba 274-8501 Japan , {\it E-mail address:watanabe.kenta@nihon-u.ac.jp}} \thanks{Department of Mathematics, Kanagawa Institute of Technology, 1030 Shimo-Ogino Atsugi city Kanagawa 243-0292 Japan, {\it E-mail address:komeda@gen.kanagawa-it.ac.jp}}}

\date{}

\maketitle 

\begin{abstract}

\noindent Let $X$ be a K3 surface, let $C$ be a smooth curve of genus $g$ on $X$, and let $A$ be a line bundle of degree $d$ on $C$. Then a line bundle $M$ on $X$ with $M\otimes\mathcal{O}_C=A$ is called a lift of $A$ . In this paper, we prove that if the dimension of the linear system $|A|$ is $r\geq2$, $g>2d-4+r(r-1)$, $d\geq 2r+4$, and $A$ computes the Clifford index of $C$, then there exists a base point free lift $M$ of $A$ such that the general member of $|M|$ is a smooth curve of genus $r$. In particular, if $|A|$ is a base point free net which defines a double covering $\pi:C\longrightarrow C_0$ of a smooth curve $C_0\subset\mathbb{P}^2$ of degree $k\geq 4$ branched at distinct $6k$ points on $C_0$, then, by using the aforementioned result, we can also show that there exists a 2:1 morphism $\tilde{\pi}:X\longrightarrow \mathbb{P}^2$ such that $\tilde{\pi}|_C=\pi$. 
 
\end{abstract}

\noindent {\bf{Keywords}} K3 extension, LM bundle, Brill-Noether theory, Donagi-Morrison lift, double covering

\smallskip

\noindent {\bf{Mathematics Subject Classification }}14J28, 14J60, 14H60

\section{Introduction}

Let $C$ be a smooth projective curve of genus $g$. Then the Brill-Noether locus $W_d^r(C)$ consisting of line bundles of degree $d$ on $C$ which have at least $r+1$ linearly independent global sections has the expected dimension $\rho(g,r,d)=g-(r+1)(g-d+r)$, where $\rho(g,r,d)$ is the Brill-Noether number. If $C$ is a general curve in the moduli space $\mathcal{M}_g$ of smooth curves of genus $g$, there is no line bundle $g_d^r$ on $C$ with $\rho(g,r,d)<0$. Therefore, a smooth curve which admits such a line bundle is called {\it{Brill-Noether special}}. In particular, if such a curve $C$ is contained in a K3 surface $X$, the problem of whether a line bundle on $C$ with negative Brill-Noether number is obtained by a restriction to $C$ of a line bundle on $X$ is interesting in the point of view of the knowledge of the relationship between the Brill-Noether theory of polarized K3 surfaces and the Brill-Noether theory of smooth curves obtained as their hyperplane sections.

Let $X$ be a K3 surface, and let $L$ be a base point free and big line bundle of sectional genus $g$ on $X$. Then the polarized K3 surface $(X,L)$ is called Brill-Noether special if there exists a non-trivial line bundle $M$ on $X$ with $M\neq L$ such that $h^0(M)h^0(L\otimes M^{\vee})\geq h^0(L)=g+1.$ If $(X,L)$ is Brill-Noether special, then any smooth curve $C\in |L|$ is Brill-Noether special. Because, the Brill-Noether number of $M\otimes\mathcal{O}_C$ is negative. However, the converse is still open.

\smallskip

\smallskip

\noindent{\bf{Conjecture 1.1}}. {\it{Let $(X,L)$ be a polarized K3 surface, and let $C\in |L|$ be a smooth irreducible curve. If $C$ is Brill-Noether special, then $(X,L)$ is Brill-Noether special.}}

\smallskip

\smallskip

\noindent The fact that if the Picard number of $X$ is one, then Conjecture 1.1 is correct is known as a classical result ([8]). Recently, it is known that Conjecture 1.1 is correct for several polarized K3 surfaces of genus $\leq 22$ (cf. [1, Theorem 1.3], [4, Theorem 1], [5, Remark 1.4]). However, in general, it is difficult to discriminate whether $(X,L)$ is Brill-Noether special for a given polarized K3 surface $(X,L)$.

For a line bundle $A$ on a smooth curve $C$ contained in a K3 surface $X$, we call a line bundle $M$ on $X$ with $M\otimes\mathcal{O}_C=A$ a {\it{lift}} of $A$. If the Brill-Noether number of $A$ is negative, and there exists a lift $M$ of $A$ satisfying $h^0(M)\geq2$, $h^0(L\otimes M^{\vee})\geq2$, and $h^1(L\otimes M^{\vee})=h^1(M)=0$, then $(X,L)$ is Brill-Noether special. However, the existence of such a lift can not always be expected. Thus, Donagi and Morrison have introduced a weak notion of a lift. If a line bundle $M$ on $X$ satisfies the following conditions, we say that $M$ is {\it{adapted}} to $|L|$.
\begin{itemize}
\item $h^0(M)\geq2$ and $h^0(L\otimes M^{\vee})\geq2$;
\item $\Cliff(M\otimes\mathcal{O}_C)$ is independent of the smooth curve $C\in |L|$.
\end{itemize} 

\noindent We note that $\deg (M\otimes\mathcal{O}_C)\leq g-1$ or $\deg (L\otimes M^{\vee}\otimes\mathcal{O}_C)\leq g-1$, and that if a line bundle $A$ on $C$ satisfies $h^0(A)\geq2$ and $\deg(A)\leq g-1$, then $A$ contributes to the Clifford index of $C$. Donagi and Morrison have conjectured the following.

\smallskip

\smallskip

\noindent{\bf{Conjecture 1.2}}. {\it{Let $(X,L)$ be a polarized K3 surface, and let $C\in |L|$ be a smooth curve of genus $g\geq2$. If a base point free line bundle $A$ on $C$ satisfying $\deg(A)\leq g-1$ has a negative Brill-Noether number, then there exists a line bundle $M$ on $X$ adapted to $|L|$ which satisfies the following conditions.
\begin{itemize}
\item $|A|$ is contained in the restriction of $|M|$ to $C$;
\item $\Cliff(M\otimes\mathcal{O}_C)\leq \Cliff(A)$.
\end{itemize} }}

\smallskip

\smallskip

\noindent The line bundle $M$ on $X$ as in Conjecture 1.2 is called a {\it{Donagi-Morrison lift}}. Lelli-Chiesa gave a necessary and sufficient condition for $A$ to have a Donagi-Morrison lift, in the case where $L$ is ample and $A$ computes the Clifford index of $C$ ([9, Theorem 1.1]). However, if $L$ is not ample, the existence of such a lift is still not mentioned. In this paper, we show that any line bundle $A$ which computes the Clifford index of a smooth curve $C\in |L|$ has a lift, under some strong conditions concerning the degree of $A$ and the genus $g$ of $C$, in the case where $L$ is not necessarily ample. Our main result is the following.

\smallskip

\newtheorem{thm}{Theorem}[section]

\begin{thm} Let $X$ be a K3 surface, let $C$ be a smooth curve of genus $g$ on $X$, and let $A$ be a complete linear system of dimension $r\geq2$ on $C$ with $g>2\deg(A)-4+r(r-1)$ and $\deg(A)\geq2r+4$. If $\Cliff(A)=\Cliff(C)$, then there exists a smooth curve $B$ of genus $r$ on $X$ with $|\mathcal{O}_C(B)|=A$. In particular, the consequence of Conjecture 1.2 is correct.\end{thm}

\smallskip

If $A$ is a base point free pencil of degree $d$ on a smooth curve $C\in |L|$ of genus $g$, and $g>\frac{d^2}{4}+d+2$, then there exists an elliptic curve $\Delta$ on $X$ such that $\mathcal{O}_X(\Delta)$ is a lift of $A$ ([10, Theorem 1]). Hence, the consequence of Conjecture 1.2 holds for $A$. On the other hand, there exists a smooth curve $C$ of genus $g=\frac{d^2}{4}+d+2$ on a K3 surface $X$ admitting a base point free pencil of degree $d$ on $C$ which does not have such a lift. Indeed, we can construct a double covering $C$ of a smooth plane curve $C_0$ of degree $k\geq4$ branched at distinct $6k$ points on $C_0$ as such a curve (Proposition 5.2). Then $C$ has a gonality pencil which admits no lift, and since $C$ is a smooth curve of genus $g=k^2+1$ and gonality $d=2k-2$, we have $g=\frac{d^2}{4}+d+2$. By Theorem 1.1, if the double covering $\pi:C\longrightarrow C_0$ is the morphism associated with a net $Z$ on $C$, there exists a lift $M$ of $Z$ defined by a smooth curve of genus 2 on $X$ (Proposition 5.4). Hence, it turns out that the morphism $\tilde{\pi}:X\longrightarrow \mathbb{P}^2$ associated with such a lift $M$ is a 2:1 map with $\tilde{\pi}|_{C}=\pi$.

Our plan of this paper is as follows. In section 2, we recall some fundamental facts concerning vector bundles on K3 surfaces. In section 3, we recall the definition of the Lazarsfeld-Mukai bundle associated with a smooth curve on a K3 surface and a base point free linear system on it, and several known results about it. In section 4, we prove Theorem 1.1. In section 5, we compute the maximum genus of smooth curves on K3 surfaces which are obtained as double coverings of smooth plane curves. Moreover, we characterize such double coverings of maximum genus.

\smallskip

\smallskip

\noindent{\bf{Notations and conventions}}. We work over the complex number field $\mathbb{C}$. In this paper, a curve and a surface are smooth and projective. Let $X$ be a curve or a surface. We denote by $K_X$ the canonical line bundle of $X$. If two divisors $D_1$ and $D_2$ on $X$ are linearly equivalent, then we will write as $D_1\sim D_2$. We denote by $|L|$ the linear system associated with a divisor or a line bundle $L$ on $X$. We call a linear system of dimension one a pencil, and call a linear system of dimension two a net. For a torsion free sheaf $E$ on $X$, we denote the rank of $E$, the dual of $E$, and the $i$-th Chern class of $E$, by $\rk(E)$, $E^{\vee}$, and $c_i(E)$, respectively.

Let $C$ be a curve. Then the gonality of $C$ is the minimum degree of pencils on $C$, and we denote it by $\gon(C)$. We note that the gonality of a smooth plane curve of degree $k\geq4$ is $k-1$. We say that a base point free line bundle $A$ on $C$ is {\it{primitive}} if $|K_C\otimes A^{\vee}|$ is base point free. For a line bundle $A$ on $C$, we denote the Clifford index of $A$ by $\Cliff(A):=\deg(A)-2\dim|A|$. We say that $A$ contributes to the Clifford index of $C$ if $h^0(A)\geq2$ and $h^1(A)\geq2$. Then the Clifford index of $C$ is the minimum value of the Clifford indices of such line bundles on $C$, and we denote it by $\Cliff(C)$. We note that $\gon(C)-3\leq\Cliff(C)\leq\gon(C)-2$. The Clifford dimension of $C$ is the minimum dimension of complete linear systems on $C$ which compute the Clifford index of $C$. If the Clifford dimension of $C$ is one, then $\Cliff(C)=\gon(C)-2$.

Let $X$ be a surface. Then we denote by $\Pic(X)$ the Picard group of $X$, and call its rank the Picard number of $X$. We denote by $\mathcal{J}_W$ the ideal sheaf of a subscheme $W$ of $X$. $X$ is called a regular surface if $h^1(\mathcal{O}_X)=0$. Moreover, if $K_X$ is trivial, then we call $X$ a K3 surface. We call a K3 surface containing a smooth curve $C$ a {\it{K3 extension}} of $C$.

\section{Vector bundles on K3 surfaces}

In this section, we recall several classical results about vector bundles on K3 surfaces. Let $X$ be a K3 surface, and let $E$ be a vector bundle on $X$. First of all, by the Riemann-Roch theorem, we obtain the following equality.

$$\chi(E)=2\rk(E)+\frac{c_1(E)^2}{2}-c_2(E),$$
where $\chi(E)=h^0(E)-h^1(E)+h^2(E)$. Since $K_X$ is trivial, by the Serre duality, we have
$$h^i(E)=h^{2-i}(E^{\vee})\; (0\leq i\leq 2).$$
In particular, if $D$ is a divisor on $X$, then
$$h^i(\mathcal{O}_X(D))=h^{2-i}(\mathcal{O}_X(-D))\; (0\leq i\leq 2).$$
Hence, if $D^2\geq-2$, then $h^0(\mathcal{O}_X(D))>0$ or $h^0(\mathcal{O}_X(-D))>0$. If $D$ is a non-zero effective divisor on $X$, then $D^2=2P_a(D)-2$, by the adjunction formula, where $P_a(D)$ is the arithmetic genus of $D$.

$\;$

\newtheorem{prop}{Proposition}[section]

\begin{prop} {\rm{([11, Proposition 2.6])}}. Let $L$ be a non-trivial line bundle on $X$. If $|L|\neq\emptyset$ and $|L|$ has no fixed component, then one of the following cases occurs. 

\smallskip

\smallskip

\noindent {\rm{(i)}} $L^2>0$ and the general member of $|L|$ is a smooth irreducible curve of genus $\frac{L^2}{2}+1$. In this case, $h^1(L)=0$.

\noindent {\rm{(ii)}} If $L^2=0$, then there exist an elliptic curve $\Delta$ on $X$ and an integer $k\geq1$ satisfying $L\cong\mathcal{O}_X(k\Delta)$. In this case, $h^1(L)=k-1$. \end{prop}

$\;$

\noindent If $C$ is an irreducible curve on $X$ with $C^2\geq0$, then $|C|$ has no base point ([11, Theorem 3.1]). Therefore, by Proposition 2.1, we obtain the following proposition.

$\;$

\begin{prop} {\rm{{([11, Corollary 3.2])}}}. Let $L$ be a non-trivial line bundle on $X$ with $|L|\neq\emptyset$. Then $|L|$ has no base point outside of its fixed components. \end{prop}

$\;$

\noindent Note that if a non-trivial line bundle $L$ on $X$ with $|L|\neq\emptyset$ has no base point, then it is nef.

\section{Lazarsfeld-Mukai bundles on K3 surfaces}

In this section, we recall the definition of a Lazarsfeld-Mukai bundle, and several important results about it. Let $X$ be a K3 surface, and let $C$ be a smooth curve of genus $g\geq2$ on $X$. Let $Z$ be a non-zero effective divisor of degree $d$ on $C$ such that $|Z|$ has no base point, and let $V\subset H^0(\mathcal{O}_C(Z))$ be a subspace which forms a base point free linear system of dimension $r\geq1$ on $C$. Then we denote by $E_{C,(Z,V)}$ the dual of the kernel of the evaluation map $V\otimes\mathcal{O}_X\longrightarrow \mathcal{O}_C(Z)$. We call it the {\it{Lazarsfeld-Mukai}} ({\it{LM}} for short) bundle on $X$ associated with $(Z,V)$. In particular, there exists the following exact sequence.
$$0\longrightarrow V^{\vee}\otimes\mathcal{O}_X\longrightarrow E_{C,(Z,V)}\longrightarrow K_C\otimes\mathcal{O}_C(-Z)\longrightarrow 0.$$
Here, we note that $V^{\vee}$ defines an $(r+1)$-dimensional subspace of $H^0(E_{C,(Z,V)})$. Since $(Z,V)$ has no base point, the evaluation map associated with it is surjective, and hence, $E_{C,(Z,V)}$ is locally free. If $V=H^0(\mathcal{O}_C(Z))$, then we denote its LM bundle merely by $E_{C,Z}$.

\smallskip

\smallskip

\begin{prop} {\rm{([9, Proposition 2.1])}}. Any LM bundle $E_{C,(Z,V)}$ on $X$ has the following properties.

\smallskip

\smallskip

\noindent {\rm{(i)}} $\rk E_{C,(Z,V)}=r+1$, $c_1(E_{C,(Z,V)})=\mathcal{O}_X(C)$, and $c_2(E_{C,(Z,V)})=d$.

\noindent {\rm{(ii)}} $h^1(E_{C,(Z,V)})=h^0(\mathcal{O}_C(Z))-r-1$, and $h^2(E_{C,(Z,V)})=0$.

\noindent {\rm{(iii)}} $E_{C,(Z,V)}$ is globally generated off the set of base points of $|K_C\otimes\mathcal{O}_C(-Z)|$.

\noindent {\rm{(iv)}} $\chi(E_{C,(Z,V)}^{\vee}\otimes E_{C,(Z,V)})=2(1-\rho(g,r,d))$. \end{prop}

\smallskip

\smallskip

\noindent If $\mathcal{O}_C(Z)$ computes the Clifford index of $C$, it is primitive. Then, by the assertion (iii) as in Proposition 3.1, $E_{C,(Z,V)}$ is generated by its global sections. Moreover, if $V=H^0(\mathcal{O}_C(Z))$, then, by the assertion (ii), we have $h^1(E_{C,Z})=h^2(E_{C,Z})=0$. We obtain the following characterization.

\smallskip

\smallskip

\begin{prop} {{\rm{([2, Lemma 1.2])}}}. If a vector bundle $E$ on $X$ satisfying $h^1(E)=h^2(E)=0$ is globally generated, then $E$ is the LM bundle $E_{C,Z}$ associated with a smooth curve $C$ on $X$ and a base point free effective divisor $Z$ on $C$. \end{prop}

\smallskip

\smallskip

\noindent By the assertion (iv) as in Proposition 3.1, if $\rho(g,r,d)<0$, then $E_{C,(Z,V)}$ is not simple. Therefore, there exists a non-zero morphism $f\in\End(E_{C,(Z,V)})$ dropping the rank everywhere. Moreover, if $r=1$ (i.e., $\rk E_{C,(Z,V)}=2$), then, by the assertion (ii) and (iii), the image of $f$ is a torsion free sheaf of rank one on $X$ whose first Chern class is not trivial and base point free. In particular, the following assertion follows.

\smallskip

\smallskip

\begin{prop} {\rm{([2, Lemma 2.1], [3, Lemma 4.4])}}. Assume that the rank of $E_{C,(Z,V)}$ is two. If $\rho(g,1,d)<0$, then there exist two line bundles $M,N\in\Pic(X)$ satisfying $h^0(M)\geq2$, $h^0(N)\geq2$, and $M^2\geq N^2$, and a zero-dimensional subscheme $W\subset X$ such that $N$ has no base point and $E_{C,(Z,V)}$ sits in the following exact sequence.
$$0\longrightarrow M \longrightarrow E_{C,(Z,V)}\longrightarrow N\otimes\mathcal{J}_W \longrightarrow0.$$
\end{prop}

\smallskip

\smallskip

\noindent By the assertion of (ii) and (iii) as in Proposition 3.1, the definition of a LM bundle on $X$ can be generalized as follows.

\smallskip

\smallskip

\newtheorem{df}{Definition}[section]

\begin{df} {\rm{([9, Definition 1])}}. Let $E$ be a torsion free sheaf on $X$ with $h^2(E)=0$. If one of the following conditions is satisfied, $E$ is called a {\rm{generalized Lazarsfeld-Mukai}} {\rm{(g.LM}} for short{\rm{) }}bundle on $X$.

\smallskip

\smallskip

\noindent {\rm{(i)}} $E$ is a locally free sheaf which is globally generated off a finite set.

\noindent {\rm{(ii)}} $E$ is globally generated. \end{df}

\smallskip

\smallskip

\noindent Next we recall the following result about a g.LM bundle $E$ on $X$ with $c_1(E)^2=0$.

\smallskip

\smallskip

\begin{prop} {\rm{([9, Proposition 2.7])}}. Let $E$ be a g.LM bundle on $X$ with $c_1(E)^2=0$. Then $E$ is a locally free sheaf with $c_2(E)=0$ generated by its global sections. Moreover, if $h^1(E)=0$, then there exists an elliptic curve $\Delta$ on $X$ satisfying $E=\mathcal{O}_X(\Delta)^{\oplus\rk E}$. \end{prop}

\smallskip

\smallskip

\noindent Since, by Proposition 3.1 (i), $\Cliff(\mathcal{O}_C(Z))=c_2(E_{C,Z})-2(\rk E_{(C,Z)}-1)$, we can generalize the notion of the Clifford index of a smooth curve on $X$ as follows, by using the notion of a g.LM bundle on $X$.

\smallskip

\smallskip

\begin{df} {\rm{([9, Definition 2])}}. We denote the Clifford index of a g.LM bundle $E$ on $X$ by $\Cliff(E):=c_2(E)-2(\rk E-1)$. \end{df}

\smallskip

\smallskip

\noindent If $E$ is a g.LM bundle on $X$, then $h^2(E^{\vee\vee})=0$ and $E^{\vee\vee}$ is globally generated off a finite set, and hence, $E^{\vee\vee}$ is also a g.LM bundle on $X$. Moreover, we obtain the following proposition.

\smallskip

\smallskip

\begin{prop} {\rm{([9, Proposition 2.4])}}. Let $E$ be a g.LM bundle on $X$ with $c_1(E)^2>0$. Then the following are satisfied:

\smallskip

\smallskip

\noindent {\rm{(i)}} If $E$ is of type {\rm{(i)}} as in Definition 3.1, then $\Cliff(E)\geq 2h^1(E)$.

\noindent {\rm{(ii)}} If $E$ is of type {\rm{(ii)}} as in Definition 3.1, then $\Cliff(E)\geq\Cliff(E^{\vee\vee})$. \end{prop}

\section{Proof of Theorem 1.1}

Let $X$ be a K3 surface, let $C$ be a smooth curve of genus $g$ on $X$, and let $A$ be a line bundle of degree $d$ on $C$ such that $\dim|A|=r\geq2$, $g>2d-4+r(r-1)$, and $d\geq2r+4$. Then we note that $\rho(g,r,d)<0$ and $d<g-1$. From now on, we assume that $A$ computes the Clifford index of $C$, and let $Z$ be an effective divisor on $C$ satisfying $A=\mathcal{O}_C(Z)$.

$\;$

\newtheorem{lem}{Lemma}[section]

\begin{lem} Let $E_{C,Z}$ be the LM bundle on $X$ associated with $C$ and $Z$. Then there exist a torsion free sheaf $F$ of rank $r$ on $X$ generated by its global sections, and a saturated line bundle $M\subset E_{C,Z}$ such that $h^0(M)\geq2$ and $E_{C,Z}$ sits in the following exact sequence.
$$0\longrightarrow M\longrightarrow E_{C,Z}\longrightarrow F \longrightarrow0.\leqno (1)$$
\end{lem}

$\;$

\noindent {\it{Proof}}. First of all, since $A$ computes the Clifford index of $C$, $A$ is base point free and primitive. By Proposition  3.1, $E_{C,Z}$ is globally generated. Here, we show that $E_{C,Z}$ contains a line bundle on $X$ which has at least two linearly independent global sections. We take a general subspace $V\in\Grass(2, H^0(A))$, where $\Grass(t, H^0(A))$ is the Grassmann manifold consisting of $t$-dimensional subspaces of $H^0(A)$. Then it forms a base point free pencil $|V|$ on $C$. We set

\smallskip

\smallskip

$F_{C,Z}:=\ker(ev:H^0(A)\otimes\mathcal{O}_X\longrightarrow A),$

$F_{C,(Z,V)}:=\ker(\tilde{ev}:V\otimes\mathcal{O}_X\longrightarrow A),$

\smallskip

\smallskip

\noindent where $ev$ and $\tilde{ev}$ are the evaluation maps associated with $H^0(A)$ and $V$, respectively. Since $|V|$ has no base point, $F_{C,(Z,V)}$ is locally free. Since $ev|_V=\tilde{ev}$, we have $F_{C,(Z,V)}\subset F_{C,Z}$. By the computation of the Chern classes of $F_{C,Z}$ and $F_{C,(Z,V)}$, we have $F_{C,Z}/F_{C,(Z,V)}\cong\mathcal{O}_X^{\oplus r-1}$. Since $E_{C,Z}=F_{C,Z}^{\vee}$ and $E_{C,(Z,V)}=F_{C,(Z,V)}^{\vee}$, we have the following exact sequence.
$$0\longrightarrow \mathcal{O}_X^{\oplus r-1}\longrightarrow E_{C,Z}\longrightarrow E_{C,(Z,V)} \longrightarrow0.\leqno (2)$$
By the hypothesis, we have $\rho(g,1,d)=2d-g-2<-r(r-1)+2$. By Proposition 3.1, we have $h^0(F_{C,(Z,V)}\otimes E_{C,(Z,V)})\geq 1-\rho(g,1,d)\geq r(r-1)$. If we apply the functor $\otimes F_{C,(Z,V)}$ to the exact sequence (2), we have the following exact sequence.
$$0\longrightarrow F_{C,(Z,V)}^{\oplus r-1}\longrightarrow F_{C,(Z,V)}\otimes E_{C,Z}\longrightarrow F_{C,(Z,V)}\otimes E_{C,(Z,V)} \longrightarrow0.$$
Since $h^1(E_{C,Z})=h^2(E_{C,Z})=0$, we get $h^1(F_{C,(Z,V)})=h^1(E_{C,(Z,V)})=r-1$, by the sequence (2). We have  $h^1(F_{C,(Z,V)}^{\oplus r-1})\leq (r-1)^2$. Since $h^0(F_{C,(Z,V)})=0$, we have $h^0(F_{C,(Z,V)}\otimes E_{C,Z})\geq h^0(F_{C,(Z,V)}\otimes E_{C,(Z,V)})-(r-1)^2\geq r-1\geq1$. Here, let $f\in \Hom(E_{C,(Z,V)}, E_{C,Z})$ be a non-zero map.

If we assume that the rank of $f$ is one, then there exist $N\in\Pic(X)$ and a zero-dimensional subscheme $W\subset X$ such that the image of $f$ coincides with $N\otimes\mathcal{J}_W$. Since, by Proposition 3.1, $E_{C,(Z,V)}$ is globally generated, $|N|$ is base point free. Moreover, since $h^2(E_{C,(Z,V)})=0$, $N$ is not trivial, and hence, $h^0(N)\geq2$. Since $(N\otimes\mathcal{J}_W)^{\vee\vee}\subset E_{C,Z}^{\vee\vee}$ and $(N\otimes\mathcal{J}_W)^{\vee\vee}=N$, we have $N\subset E_{C,Z}$.

Assume that the rank of $f$ is two. Since $\rho(g,1,d)<0$, by Proposition 3.3, there exist line bundles $M,N\in\Pic(X)$ and a zero-dimensional subscheme $W\subset X$ such that $h^0(M)\geq2$, $h^0(N)\geq2$, and $E_{C,(Z,V)}$ sits in the following exact sequence.
$$0\longrightarrow M\longrightarrow E_{C,(Z,V)}\longrightarrow N\otimes\mathcal{J}_W \longrightarrow0.$$
Since $f$ is injective, we have $M\subset E_{C,Z}$. By the above observation, there exists a saturated line bundle $M\subset E_{C,Z}$ satisfying $h^0(M)\geq2$. Since $E_{C,Z}$ is globally generated, if we set $F=E_{C,Z}/M$, then $F$ is a torsion free sheaf of rank $r$ on $X$ generated by its global sections. $\hfill\square$

$\;$

\begin{lem} Notations be as in Lemma 4.1. Then there exist a smooth curve $B$ on $X$, and an effective divisor $Z^{'}$ of degree $2r-2$ on $B$ which forms a base point free linear system of dimension $r-1$ on $B$ such that $F\cong E_{B,Z^{'}}$. Moreover, $h^1(M)=0$. \end{lem}

$\;$

\noindent {\it{Proof}}. The torsion free sheaf $F$ on $X$ as in Lemma 4.1 sits in the following exact sequence.
$$0\longrightarrow F\longrightarrow F^{\vee\vee}\longrightarrow S \longrightarrow0,\leqno (3)$$
where $S$ is a coherent sheaf of finite length on $X$. We set $c_1(F^{\vee\vee})=N$ and let $W$ be the support of $S$. Then, by [6, (0.3)], we have
$$c_2(E_{C,Z})=\deg Z=d=M.N+|W|+c_2(F^{\vee\vee}),\leqno (4)$$
where $|W|$ is the length of $W$. Since $h^2(E_{C,Z})=0$, by the exact sequence (1), we have $h^2(F)=0$. Moreover, by the sequence (3), we have $h^2(F^{\vee\vee})=0$. On the other hand, since $h^1(E_{C,Z})=0$, by the exact sequence (1), we have $h^1(F)=0$. Since $F$ is generated by its global sections, $N$ is base point free. Indeed, since $c_1(F^{\vee\vee})=\det(F^{\vee\vee})$ is base point free off the set $W$, by Proposition 2.2, it is base point free everywhere.

Assume that $N^2=0$. By Proposition 3.4, there exists an elliptic curve $\Delta$ on $X$ such that $F=F^{\vee\vee}\cong\mathcal{O}_X(\Delta)^{\oplus r}$. Then we have $c_2(F^{\vee\vee})=0$ and $|W|=0$. Since $N\cong\mathcal{O}_X(r\Delta)$, $\deg (\mathcal{O}_C(r\Delta))=M.N=d$. Since $h^0(\mathcal{O}_C(\Delta))\geq2$ and $\gon(C)\geq d-2r+2$, we have $d-2r+2\leq\deg(\mathcal{O}_C(\Delta))=\frac{d}{r}$. However, this contradicts the assumption $d\geq 2r+4$. Hence, we have
$$N^2>0.\leqno (5)$$

On the other hand, since $h^0(M)\geq2$ and $h^0(N)\geq2$, $M\otimes\mathcal{O}_C$ contributes to the Clifford index of $C$. Since $N$ has no base point, by the inequality (5) and Proposition 2.1, we have $h^1(N)=0$, and hence, we have $h^0(N\otimes\mathcal{O}_C)=h^0(N)+h^1(M)$. Therefore, by the equality (4), we have

\smallskip

\smallskip

\noindent $d-2r=\Cliff(C)\leq\Cliff(M\otimes\mathcal{O}_C)=\Cliff(K_C\otimes M^{\vee})=\Cliff(N\otimes\mathcal{O}_C)$

$=N.C-2(h^0(N)+h^1(M)-1)$

$=N.C-N^2-2-2h^1(M)=M.N-2-2h^1(M)$

$=d-|W|-c_2(F^{\vee\vee})-2-2h^1(M).$

\smallskip

\smallskip

\noindent Hence, we have
$$|W|+c_2(F^{\vee\vee})+2h^1(M)\leq2r-2.\leqno (6)$$
Since $h^2(F^{\vee\vee})=0$ and $F^{\vee\vee}$ is globally generated off a finite set, by the inequality (5) and Proposition 3.5, we have $c_2(F^{\vee\vee})\geq2r-2$. By the inequality (6), we obtain $|W|=0$, $h^1(M)=0$, $c_2(F^{\vee\vee})=2r-2$ and $h^1(F^{\vee\vee})=0$. In particular, since $F=F^{\vee\vee}$, $F^{\vee\vee}$ is globally generated and $h^1(F^{\vee\vee})=0$. By Proposition 3.2, there exist a smooth curve $B\in |N|$ and an effective divisor $Z^{'}$ of degree $2r-2$ on $B$ such that $|Z^{'}|$ is a base point free linear system of dimension $r-1$ on $B$, and $F^{\vee\vee}\cong E_{B,Z^{'}}$. $\hfill\square$

\smallskip

\smallskip

\begin{lem} Let the notations be as in Lemma 4.2. Then the genus of $B$ is $r$ or $r-1$. \end{lem}

\smallskip

\smallskip

\noindent{\it{Proof}}. Let $g(B)$ be the genus of $B$, and let $Z^{'}$ be as in Lemma 4.2. Since $h^0(\mathcal{O}_B(Z^{'}))=r$, by the Riemann-Roch theorem, we have $g(B)\geq r-1$. Assume that $g(B)>r$. Since $h^1(\mathcal{O}_B(Z^{'}))\geq2$, $\mathcal{O}_B(Z^{'})$ contributes to the Clifford index of $B$. Since $\Cliff(\mathcal{O}_B(Z^{'}))=0$, we have $\Cliff(B)=0$. By [7, Theorem 2, Proposition 4.1], $B$ is a hyperelliptic curve. Since $r\geq2$, we have $g(B)\geq3$. By [11, Theorem 5.2], one of the following cases occurs.

\smallskip

\smallskip

\noindent (i) There exists a smooth genus 2 curve $B_0$ on $X$ such that $B\sim 2B_0$.

\noindent (ii) There exists an elliptic curve $\Delta$ on $X$ such that $B.\Delta=2$.

\smallskip

\smallskip

We consider the case of (i). First of all, we have $B^2=8$. Since $c_1(F)=\mathcal{O}_X(B)$, we obtain $C\sim M+B$. By Lemma 4.1, Lemma 4.2, and the equality (4), we have $M.B=d-2r+2$. Hence,  $C.B=M.B+B^2=d-2r+10$. Since $h^0(\mathcal{O}_C(B_0))\geq3$ and 
$$h^0(K_C\otimes\mathcal{O}_C(-B_0))=h^0(\mathcal{O}_X(C-B_0))>h^0(M)\geq2,$$
$\mathcal{O}_C(B_0)$ contributes to the Clifford index of $C$. Since $\deg(\mathcal{O}_C(B_0))=\frac{d}{2}-r+5$, we have $d-2r\leq\Cliff(\mathcal{O}_C(B_0))\leq\frac{d}{2}-r+1$. Therefore, we have $d\leq2r+2$. However, this contradicts the hypothesis $d\geq2r+4$.

\smallskip

\smallskip

We consider the case of (ii). Since $B^2$ is even and $B^2\geq4$, there exists a natural number $s\in\mathbb{N}$ such that $B^2=4s$ or $4s+2$. Then we have $(B-s\Delta)^2=B^2-4s\geq0$ and $\Delta.(B-s\Delta)=2>0$. Since $\Delta$ is nef, we have $h^0(\mathcal{O}_X(B-s\Delta))\geq2$. By the exact sequence
$$0\longrightarrow \mathcal{O}_X(-s\Delta)\otimes M^{\vee}\longrightarrow \mathcal{O}_X(B-s\Delta)\longrightarrow \mathcal{O}_C(B-s\Delta) \longrightarrow0,$$
we have $h^0(\mathcal{O}_C(B-s\Delta))\geq h^0(\mathcal{O}_X(B-s\Delta))\geq \chi(\mathcal{O}_X(B-s\Delta))=\frac{B^2}{2}-2s+2$. Since $\gon(C)\geq d-2r+2$, we have
$$C.(B-s\Delta)-(d-2r+2)\geq h^0(\mathcal{O}_C(B-s\Delta))-2\geq\frac{B^2}{2}-2s.$$
Since $C.B=B^2+B.M=B^2+d-2r+2$ and $B.\Delta=2$, we obtain $B^2\geq 2sM.\Delta$. Hence, we have
$$M.\Delta\leq 3.\leqno (7)$$
On the other hand, since $\Delta.(\Delta-C)=-\Delta.C<0$, we have $h^0(\mathcal{O}_X(\Delta-C))=0$, and hence, by the exact sequence
$$0\longrightarrow \mathcal{O}_X(\Delta-C)\longrightarrow \mathcal{O}_X(\Delta)\longrightarrow \mathcal{O}_C(\Delta) \longrightarrow0,$$
we have $h^0(\mathcal{O}_C(\Delta))\geq h^0(\mathcal{O}_X(\Delta))=2$. Since $\gon(C)\geq d-2r+2$, we have $C.\Delta\geq d-2r+2$, and hence, $M.\Delta\geq d-2r$. By the inequality (7), we have $d\leq 2r+3$. This contradicts the hypothesis $d\geq 2r+4$. By the above observation, we have $g(B)=r$ or $r-1$. $\hfill\square$

\smallskip

\smallskip

{\it{Proof of Theorem 1.1}}. If we apply the functor $\otimes M^{\vee}$ to the exact sequence
$$0\longrightarrow H^0(A)^{\vee}\otimes\mathcal{O}_X\longrightarrow E_{C,Z}\longrightarrow K_C\otimes A^{\vee}\longrightarrow0,$$
we obtain the exact sequence
$$0\longrightarrow H^0(A)^{\vee}\otimes M^{\vee}\longrightarrow E_{C,Z}\otimes M^{\vee}\longrightarrow N\otimes\mathcal{O}_C\otimes A^{\vee}\longrightarrow0.$$
Since, by Lemma 4.2, $h^1(M)=0$, we have $h^0(N\otimes\mathcal{O}_C\otimes A^{\vee})=h^0(E_{C,Z}\otimes M^{\vee})$. On the other hand, since $M\subset E_{C,Z}$, we have $h^0(E_{C,Z}\otimes M^{\vee})>0$. In particular, we have $N.C\geq\deg(A)$. By Lemma 4.2, the first Chern class $N$ of the sheaf $F$ which appears in the exact sequence (1) coincides with $\mathcal{O}_X(B)$. Since, by Lemma 4.3, $B^2=2r-2$ or $2r-4$, we obtain
$$\deg(A)=M.N+2r-2\geq M.N+N^2=N.C=\deg(\mathcal{O}_C(B)).$$
Therefore, we have $N\otimes\mathcal{O}_C\cong A$, and the genus of $B$ is $r$. $\hfill\square$

$\;$

\noindent{\bf{Corollary 4.1}}. {\it{Let $X$ be a K3 surface, and let $C$ be a smooth curve of genus $g$ on $X$. We set $L=\mathcal{O}_X(C)$. If there exists a complete linear system $A$ of degree $d$ and dimension $r\geq2$ on $C$ which computes the Clifford index of $C$ such that $g>2d-4+r(r-1)$ and $d\geq2r+4$, then the polarized K3 surface $(X,L)$ is Brill-Noether special.}}

$\;$

\noindent {\it{Proof}}. By Theorem 1.1, there exists a base point free line bundle $B$ of sectional genus $r$ on $X$ satisfying $|B\otimes\mathcal{O}_C|=A$. Since $B$ is nef and $d\geq 2r+4$, by the Serre duality, we have $h^2(L\otimes B^{\vee})=0$. By the Riemann-Roch theorem, we have $h^0(L\otimes B^{\vee})\geq g-d+r$. Hence, we have 
$$h^0(B)h^0(L\otimes B^{\vee})-h^0(L)\geq (r+1)(g-d+r)-g-1=-\rho(g,r,d)-1>0.$$
Therefore, we obtain the consequence. $\hfill\square$

\section{Double coverings of plane curves admitting K3 extensions}

In this section, we compute the maximum genus of double coverings of smooth plane curves contained in K3 surfaces. Moreover, we characterize such double coverings of maximum genus. First of all, we prove the following proposition.

$\;$

\begin{prop} Let $\pi:C\longrightarrow C_0$ be a double covering of a smooth plane curve $C_0$ of degree $k\geq3$, and let $g$ be the genus of $C$. If $C$ is contained in a K3 surface, then $g\leq k^2+1$. \end{prop}

\smallskip

\smallskip

\noindent {\it{Proof}}. Let $d$ be the gonality of $C$, and let $g_0$ be the genus of $C_0$. Since the gonality of $C_0$ is $k-1$, we have $d\leq 2k-2$. Assume that $g>k^2+1$. Since $g_0=\frac{1}{2}(k-1)(k-2)$, we have $d<g-2g_0$. Hence, by the inequality of Castelnuovo-Severi, we have $d=2k-2$. Since $k^2+1=\frac{d^2}{4}+d+2$, any gonality pencil on $C$ has a lift, by [10, Theorem 1]. However, since $\dim W_d^1(C)>0$ and any K3 surface has at most countably many elliptic pencils, this is a contradiction. $\hfill\square$

$\;$

\begin{prop} Let $C_0$ be a smooth plane curve of degree $k\geq3$. Then there exists a double covering $C$ of genus $g=k^2+1$ of $C_0$ which has a K3 extension. \end{prop}

\smallskip

\smallskip

\noindent {\it{Proof}}. Let $D\subset\mathbb{P}^2$ be a smooth curve of degree six which intersects transversely with $C_0$ at distinct $6k$ points, and let $\tilde{\pi}:X\longrightarrow\mathbb{P}^2$ be a double covering branched along $D$. Moreover, let $C$ be the fiber product of $\iota : C_0\hookrightarrow\mathbb{P}^2$ and $\tilde{\pi}$. Then $C$ is a double covering of $C_0$ branched at $D\cap C_0$. By the Hurwitz formula, the genus of $C$ is $k^2+1$. $\hfill\square$

$\;$

\noindent By Proposition 5.2, the maximum genus $g$ of double coverings of smooth plane curves which have K3 extensions is $\frac{d^2}{4}+d+2$, where $d$ is the gonality of such a double covering which has the maximum genus $g$. Conversely, any smooth curve $C$ of genus $g=\frac{d^2}{4}+d+2$ contained in a K3 surface admitting a base point free pencil of degree $d\geq4$ which has no lift can be characterized as follows.

$\;$

\begin{prop} Let $X$ be a K3 surface, and let $C$ be a smooth curve of genus $g$ on $X$. Assume that there exists a base point free divisor $Z$ of degree $d\geq4$ on $C$ such that $\mathcal{O}_C(Z)$ is a pencil on $C$ which has no lift, and $g=\frac{d^2}{4}+d+2$. Then there exists a smooth genus 2 curve $B$ on $X$ with $C\sim kB$, where $k$ is an integer with $d=2k-2$. \end{prop}

\smallskip

\smallskip

\noindent {\it{Proof}}. Since $\rho(g,1,d)<0$, the LM bundle $E_{C,Z}$ on $X$ associated with $C$ and $Z$ is not simple. By Proposition 3.3, there exist line bundles $M,N\in \Pic(X)$ satisfying $M^2\geq N^2$, $h^0(M)\geq2$, and $h^0(N)\geq2$, and a zero-dimensional subscheme $W\subset X$ such that $N$ is base point free and $E_{C,Z}$ sits in the following exact sequence.
$$0\longrightarrow M \longrightarrow E_{C,Z}\longrightarrow N\otimes\mathcal{J}_W \longrightarrow0.\leqno (8)$$

\smallskip

\smallskip

\begin{lem} The line bundle $N$ on $X$ which appears in the exact sequence {\rm{(8)}} satisfies $N^2>0$.\end{lem}

\smallskip

\smallskip

\noindent {\it{Proof}}. Since $h^0(M^{\vee})=0$, if we apply the functor $\otimes M^{\vee}$ to the exact sequence
$$0\longrightarrow H^0(\mathcal{O}_{C}(Z))^{\vee}\otimes\mathcal{O}_X\longrightarrow E_{C,Z}\longrightarrow K_{C}\otimes\mathcal{O}_{C}(-Z)\longrightarrow0,$$
we have $h^0(N\otimes\mathcal{O}_{C}(-Z))\geq h^0(E_{C,Z}\otimes M^{\vee})>0$. Hence, $C.N\geq d=\deg(Z)$. 

By the computation of the Chern classes of $E_{C,Z}$, we have $\mathcal{O}_X(C)\cong M\otimes N$ and $M.N\leq d$. If $N^2=0$, then $C.N\leq d$. Hence, we have $N\otimes\mathcal{O}_{C}\cong\mathcal{O}_{C}(Z)$. This contradicts the assumption of Proposition 5.3. $\hfill\square$

$\;$

\noindent Since $M^2\geq N^2$, by Lemma 5.1, we have $M^2>0$. By the Hodge index theorem, we have $M^2N^2\leq d^2$. Since $N^2$ is even, we have $2\leq N^2\leq d$. 

Here, we set $s=N^2$ and $f(s)=(s+d)^2-C^2s$. By the hypothesis concerning the genus of $C$, we have $C^2=\frac{d^2}{2}+2d+2$, and hence, $f(s)=s^2-(\frac{d^2}{2}+2)s+d^2$. Since $d\geq4$, we have $\frac{d^2}{4}+1>d$. Hence, we have $f(s)\leq0$, and the equality holds only if $s=2$.

On the other hand, since, by the Hodge index theorem, $C^2s\leq (N.C)^2\leq(s+d)^2$, we have $f(s)\geq0$. Hence, we obtain $f(s)=0$ and $s=2$. Since $N$ is base point free, by the theorem of Bertini, the general member of $|N|$ is a smooth genus 2 curve. Since $g=\frac{d^2}{4}+d+2$, $d$ is even. Hence, there exists an integer $k\geq3$ such that $d=2k-2$. If we let $B\in |N|$ be a smooth genus 2 curve on $X$, then we have $(C-kB)^2=0$. By the Riemann-Roch theorem, we have $h^0(\mathcal{O}_X(C-kB))>0$ or $h^0(\mathcal{O}_X(kB-C))>0$. Since $C.(C-kB)=0$, if the first case occurs, then by the exact sequence
$$0\longrightarrow \mathcal{O}_X(-kB)\longrightarrow \mathcal{O}_X(C-kB)\longrightarrow \mathcal{O}_{C} \longrightarrow0,$$
we have $h^0(\mathcal{O}_X(C-kB))=1$. Since $B.(kB-C)=0$, if the latter case occurs, then by the exact sequence
$$0\longrightarrow \mathcal{O}_X((k-1)B-C)\longrightarrow \mathcal{O}_X(kB-C)\longrightarrow \mathcal{O}_B \longrightarrow0,$$
we have $h^0(\mathcal{O}_X(kB-C))=1$. By the above observation, we get $C\sim kB$. $\hfill\square$

$\;$

\noindent{\bf{Corollary 5.1}}. {\it{Let $C$ be a double covering of a smooth plane curve $C_0$ of degree $k\geq3$ branched at distinct $6k$ points on $C_0$, and assume that $C$ is contained in a K3 surface $X$. Then there exists a smooth genus 2 curve $B$ on $X$ with $C\sim kB$.}}

$\;$

\noindent {\it{Proof}}. The gonality of $C$ is $2k-2$ and the genus of $C$ is $k^2+1$. By the same reason as in the proof of Proposition 5.1, there exists a gonality pencil on $C$ which has no lift. Hence, by Proposition 5.3, the assertion holds. $\hfill\square$

$\;$

\begin{prop} Let $C_0$ be a smooth plane curve of degree $k\geq4$, and let $\pi:C\longrightarrow C_0$ be a double covering branched at distinct $6k$ points on $C_0$ such that $C$ is contained in a K3 surface $X$. If $\pi$ is the morphism associated with a base point free net $A$ of degree $2k$ on $C$, then there exists a smooth genus 2 curve $B$ on $X$ such that $\mathcal{O}_X(B)$ is a lift of $A$. \end{prop}

$\;$

\noindent In Proposition 5.4,  since the genus of $B$ is two,  $h^0(A)=h^0(\mathcal{O}_X(B))$. Hence, the consequence means that if we set $\tilde{\pi}=\Phi_{|B|}:X\longrightarrow\mathbb{P}^2$, then $\tilde{\pi}|_{C}=\pi$. From now on, let $C$ and $A$ be as in Proposition 5.4.

$\;$

\begin{lem} The Clifford dimension of $C$ is one. \end{lem}

\smallskip

\smallskip

\noindent {\it{Proof}}. Since $k\geq4$, the linear system $|C|$ on $X$ is not the Donagi-Morrison's example  (see [2, Theorem A] or  [7]). Assume that there exist a smooth curve $D$ of genus $\geq2$ and a smooth rational curve $\Gamma$ on $X$ satisfying $C\sim 2D+\Gamma$ and $D.\Gamma=1$. Since $C.\Gamma=0$, we have $4D^2=(C-\Gamma)^2={C}^2-2$. By Corollary 5.1, there exists a smooth genus 2 curve $B_0$ on $X$ with $C\sim kB_0$. Therefore, we have $D^2=\frac{1}{2}(k^2-1)$. Since $k\geq4$, we have $B_0^2=2<D^2-1$ and $0<C.B_0-B_0^2=2k-2<D^2$. Hence, by [7, Theorem 2], we obtain the consequence of Lemma 5.2. $\hfill\square$

\smallskip

\smallskip

\noindent {\it{Proof of Proposition 5.4}}. By Lemma 5.2, $\Cliff(A)=2k-4=\gon(C)-2=\Cliff(C)$. Since the genus of $C$ is $k^2+1$ and $k\geq 4$, there exists a smooth genus 2 curve $B\subset X$ such that $\mathcal{O}_X(B)$ is a lift of $A$, by Theorem 1.1. $\hfill\square$

\end{document}